\documentclass[12pt, twoside]{amsart}
\usepackage{amsfonts,graphicx,amssymb}
\usepackage{amsthm}

\newtheorem{theorem}{Theorem}

\newtheorem*{definition}{Definition}

\def\bbr{{\mathbb R}}

\def\bbc{{\mathbb C}}

\def\part{\partial}

\begin{document}

\title[Comparing volumes by concurrent cross-sections of complex lines]
{Comparing volumes by concurrent cross-sections of complex lines: a Busemann-Petty type problem}

\author{Eric L.~Grinberg} 
\address{
University of Massachusetts, Boston \\
Boston MA 02125, USA}
\email{eric.grinberg@umb.edu}


\maketitle

\begin{abstract}
We consider the problem of comparing the volumes of two star bodies in an even-dimensional euclidean space $\mathbb R^{2n} = \mathbb C^n$ by comparing their cross sectional areas along complex lines (special 2-dimensional real planes) through the origin. Under mild symmetry conditions on one of the bodies a Busemann-Petty type theorem holds.  Quaternionic and Octonionic analogs also hold. The argument  relies on integration in polar coordinates coupled with Jensen's inequality. Along the way we provide a criterion that detects which centered bodies are {\it circular}. i.e., stabilized by multiplication by complex numbers of unit  modulus. Our goal is to present a Busemann-Petty type result with a minimum of required background (in the spirit of L.K. Hua's book on the classical domains) and, in addition, to suggest characterizations of classes of star bodies by means of integral geometric inequalities.
\end{abstract}

\bigskip

\par

\noindent
{
\it Domains with mild symmetry \hfil}

\bigskip

\noindent
\begin{definition} Let $D$ be domain in $\bbc^n$ and assume tacitly that all vectors under arrows belong to $\mathbb C^n$.
We say that $D$ is {\bf a circular domain, with center $\, \overrightarrow 0 \, $} if the following holds:

\[
\overrightarrow z \in D , \, \theta \in \bbr  
\quad \Longrightarrow \quad
         ( e^{\sqrt{-1}\theta}) \overrightarrow z \in D.
\]
$D$ is said to be {\it circular} with {\it center} $\overrightarrow c$
when 
$D-\overrightarrow c$ is circular with center 
$\overrightarrow 0$. 

\noindent

\end{definition}

\noindent
Thus the origin centered unit ball $\{ \vec z \, \Big\vert  \| \vec z \| = 1 \}$
and the standard unit polydisc,  
$$
P \equiv \{ \vec z  \,\,  \Big\vert \, \,  \vert  z_i \vert \le 1, i=1 \ldots n \},
$$ 
are both circular.

In the convex geometry literature circular bodies are sometimes said to have $R_\theta$ symmetry. In the functional analysis literature one sometimes speaks of a {\it balanced
set}.

\begin{theorem}  Let $D$ be a star-shaped origin symmetric
domain in $\bbc^n$, let $G_{1,n}$ be the set of $1$-dimensional vector subspaces of $\mathbb C^n$, that is, the (Grassmannian)  set of complex lines through $\overrightarrow 0$, and let $d \ell$ be the standard probability measure on $G_{1,n}$.  Then
\begin{enumerate}
\item 
\(
\textrm{vol}(D) \ge \frac{1}{n!}  \int_{ \ell \in G_{1,n}} 
	\left( \textrm{area}(D \cap \ell) \right)^n \,d\ell.
\)
\item
\(
\textrm{vol}(D) = \frac{1}{n!}  \int_{ \ell \in G_{1,n}} 
	\left( \textrm{area}(D \cap \ell) \right)^n \,d\ell
\text{ precisely when $D$ is circular.}
\)
\end{enumerate}
\end{theorem}
\noindent
Loosely we say that the volume of $D$ dominates the  mean of the $n^{th}$ power of cross-sectional areas of $D$ by lines in $G_{1,n}$, with equality precisely when $D$ is circular.

This suggests the quest of characterizing all `classes' of bodies by  geometric integral inequalities and identities, where {\it classes}, {\it geometric}, and {\it integral} are left to be made more precise (not to mention {\it all }). For instance,  characterizations of centered bodies in $\mathbb R^n$ are found in \cite{Gardner} and  \cite{Schneider}. We seek a mapping (functor?)  from classes of bodies to geometric inequalities, perhaps in the spirit of  \cite{Zalc}, where a parallel is drawn between mean value properties and differential equations. For some results in this direction see \cite{Dann-Zymonopoulou}. As an intermediate goal, one may aim to characterize the {\it bounded symmetric hermitian domains} or {\it Siegel domains} in the style and spirit of L.K.~Hua's book \cite{Hua}.

It would also be interesting to provide \emph{stability estimates}. That is, if the inequality in Theorem 1 is nearly an equality is the body in question nearly circular, perhaps as measured by some kind of \emph{modulus of circularity}? For extensive discussions of such estimates in convex geometry and geometric inequalities see \cite{Groemer} and \cite{Toth}.

\begin{proof}
Viewing the ambient space $\bbc^n$ as $\bbr^{2n}$, with its unit sphere $S^{2n-1}$, we can compute the volume of the region $D$ in polar coordinates. We recall that \( \rho_D ( \omega ) : S^{2n-1} \longrightarrow \mathbb R $, the {\it radial function} of $D$,   is defined by 
$$
\rho_D ( \omega ) = \max \{ t \in [0, + \infty ) \, \vert \, t \omega \in D \}.
$$
With this notation we have

$$
\textrm{vol}(D) = \frac{1}{2n} \int_{  \omega \in S^{2n-1}} \left(\rho_D( \omega )\right)^{2n} 
			\, d \omega ,
$$
where $ d \omega$ is  the usual surface measure on the sphere $S^{2n-1}$  induced by the Euclidean metric. Since multiplication by a phase factor stabilizes the sphere and preserves its measure, 
$$
\int_{  \omega \in S^{2n-1}} f(\omega )	\, d \omega = 
\int_{  \omega \in S^{2n-1}} f( e^{i \theta}\omega )	\, d \omega
$$
for any $\theta$.

Let \( d \theta \) be the probability measure on the unit circle. Applying Jensen's inequality \cite{Rudin} to the convex function $x \mapsto x^n$ on $[0, +\infty )$ we obtain
\[
  \int_{\theta \in S^1} \left\{ \left(\rho_D( e^{ i \theta } \omega )\right)^{2}
		               \right\}^n \,  d\theta 
\, \, \ge \,\,
  \left\{ \int_{\theta \in S^1}  \left(\rho_D( e^{ i \theta } \omega )\right)^{2}
   \,  d\theta \right\}^n.
\]

With these observations we have
$$
\begin{array}{lcl}
\textrm{vol}(D) &=& \frac{1}{2n}  \int_{  \omega \in S^{2n-1}} \left(\rho_D( \omega )\right)^{2n} 
											\, d \omega \\ \\
		&=& \frac{1}{2n} \int_{  \omega \in S^{2n-1}} 
		                 \left(\rho_D( e^{ i \theta } \omega )\right)^{2n} 
		               										\, d \omega \\ \\
                 &=& \frac{1}{2n}  \int_{  \omega \in S^{2n-1}} 
		              \  \int_{\theta \in S^1} \left(\rho_D( e^{ i \theta } \omega )\right)^{2n} 										                                                                               \,  d\theta  \, d \omega \\
										\\
		&=& \frac{1}{2n} \int_{  \omega \in S^{2n-1}} 
		               \int_{\theta \in S^1} \left\{ \left(\rho_D( e^{ i \theta } \omega )\right)^{2}
		               \right\}^n
		                										\, 
		                  d\theta 
										\, d \omega \\
										\\
		&\ge & \frac{1}{2n} \int_{  \omega \in S^{2n-1}} 
		               \left\{ \int_{\theta \in S^1}  \left(\rho_D( e^{ i \theta } \omega )\right)^{2}
		                \, 
		                d\theta \right\}^n 
										\, d \omega \\ 
										\\
		&=& \frac{1}{2n} \int_{  \omega \in S^{2n-1}} 
		               \left\{ \frac{1}{\pi} area ( D \cap \ell (\omega )) \right\}^n 
										\, d \omega \, ,\\
										
\end{array}			
$$
where the inequality is Jensen's.  Notice that equality holds iff $\rho_D( e^{ i \theta } \omega ) = \rho_D(  \omega )$ for almost all (and hence for all) $\theta$ and $\omega$ (by continuity), i.e., if and only if  $D$ is a circular region. 

This nearly completes the proof, but we stated the Theorem in terms of integration over the Grassmannian $G_{1,n}$ (with  probability measure) and we estimated volume by integrating over the sphere $S^{2n-1}$ , the Stiefel manifold if you will (with surface area measure). The two integrals are related by a constant that turns out to be  \(  2n ( \pi^n  /n!) \), \cite[p. \emph{xxi}]{Schneider}. Thus
\[
\begin{array}{lcl}
\textrm{vol}(D) & \ge &
\frac{1}{2n} \int_{  \omega \in S^{2n-1}} 
		               \left\{ \frac{1}{\pi} area ( D \cap \ell (\omega )) \right\}^n 
										\, d \omega  \\
&=&
\frac{1}{2n} \cdot 2n \cdot \dfrac{\pi^n}{n!}   \int_{  \ell \in G_{1,n}} 
		               \left\{ \frac{1}{\pi} area ( D \cap \ell ) \right\}^n 
										\, d \omega 
= 		
\dfrac{1}{n!}  \int_{ \ell \in G_{1,n}} 
	\left( \textrm{area}(D \cap \ell) \right)^n \,d\ell	.
\end{array}				
\]
\end{proof}
\noindent

\begin{theorem} 
Let $A$ and $B$ be origin-symmetric star
bodies in $\bbc^n$  with $A$ circular. If for every complex line $\ell$
through the origin in $\bbc^n$ we have 
$$
\textrm{area}(A \cap \ell) 
\, \le  \, \textrm{area} (B \cap \ell),
$$
\newline
then $A$ has smaller or equal volume compared to  $B$.
\end{theorem}

\medskip

\noindent
This is a variation on the classical {\it Busemann-Petty} problem.
It just involves integration over polar coordinates (adapted to complex geometry) together with Jensen's inequality. The general Busemann-Petty problem has a vast literature and its solutions are powered by the notion of \emph{intersection body}, introduced by E.~Lutwak in the paper \cite{Lutwak}, which has played a decisive role in the subject, and for which this result may be viewed as a small manifestation.

\begin{proof} By the previous theorem (statement and proof),

\[
\begin{array}{lcl}  \\
\textrm{vol}(A) & = & \frac{1}{n!} \int_{G_{1,n} }
	\left( \textrm{area}(A \cap \ell) \right)^n \,d\ell
	\\ \\
	              &  \le  & \frac{1}{n!}  \int_{G_{1,n}}
	\left( \textrm{area}(B \cap \ell) \right)^n \,d\ell 
	\\ \\
	 & \le & \textrm{vol}(B).            
\end{array}
\]
\end{proof}

\medskip
\noindent
Note: a similar argument allows comparison of bodies $K,L$ in $\mathbb R^{4n}$ by cross-sections along quaterionic lines (assuming $K$ is quaternion-circular, or $S^3$-circular.)

\section*{Acknowledgement}

The author wishes to thank Susanna Dann, David Feldman, Daniel Klain, Erwin Lutwak, Mehmet Orhon, Larry Zalcman and others for helpful discussions and the referees for suggesting several improvements in the manuscript.

\end{document}